\documentclass[a4paper]{article}

\usepackage{amsmath}
\usepackage{amssymb}
\usepackage{color}
\usepackage[dvips]{graphics}
\usepackage{graphicx}
\usepackage{cite}

\newtheorem{property}{Property}

\definecolor{light}{gray}{.85}

\title{Pad\'{e}--type rational and barycentric interpolation}

\author{Claude Brezinski\thanks{Laboratoire Paul Painlev\'e, UMR CNRS 8524, UFR
de Math\'ematiques Pures et Appliqu\'ees, Universit\'e des Sciences
et Technologies de Lille,
59655--Villeneuve d'Ascq cedex,
France, E--mail: {\tt Claude.Brezinski@univ-lille1.fr}.}
\and Michela Redivo--Zaglia\thanks{Universit\`a degli Studi di Padova,
Dipartimento di Matematica,
Via Trieste 63, 35121--Padova,
Italy. E--mail: {\tt Michela.RedivoZaglia@unipd.it}.}
}

\begin{document}

\maketitle

\begin{abstract}
In this paper, we consider the particular case of the general rational Hermite interpolation problem where
only the value of the function is interpolated at some points, and where the function and its first derivatives
agree at the origin. Thus, the interpolants constructed in this way possess a Pad\'{e}--type property at 0.
Numerical examples show the interest of the procedure. The interpolation procedure can be easily
modified to introduce a partial knowledge on the poles and the zeros of the function to approximated.
A strategy for removing the spurious poles is explained.
A formula for the error is proved in the real case. Applications are given.
\end{abstract}



\section{Pad\'{e}--type approximation and rational interpolation}

For representing a function $f$, rational functions are usually more powerful than polynomials.
The information on the function $f$ can consist either in the first coefficients of its Taylor series expansion around zero,
or in its values at some points of the complex plane.

\vskip 2mm

In the first case, Pad\'{e}--type, Pad\'{e}, or partial Pad\'{e} approximants can be used. They are rational functions
whose series expansion around zero (obtained by Euclidean division in ascending powers of the numerator by the denominator)
coincides with the series $f$ as far as possible. In Pad\'{e}--type approximation, the denominator can be arbitrarily chosen and, then,
the coefficients of the numerator are obtained by imposing the preceding approximation--through--order conditions. In Pad\'{e}
approximation, both the denominator and the numerator are fully determined by these conditions. For partial
Pad\'{e} approximants, a part of the denominator and/or a part of the numerator can be arbitrarily chosen, and their remaining parts are
given by the approximation--through--order conditions. On these topics, see \cite{birk,hand,iser,bak}.

\vskip 2mm

In the second case, an interpolating rational function can be built using Thiele's formula, which comes out
from continued fractions (see, for example, \cite[pp. 102ff.]{cbmrz} or \cite[Sec. III.3--4]{cuyt}).
It achieves the maximum number of interpolation conditions, and, so, no choice is left for its construction \cite{cbmrz}.
The same is true for Hermite rational interpolants, a subject treated in many publications (see, for example, \cite{csww2})
which is related to Newton--Pad\'{e} approximants \cite[p. 157]{cuyt}.
On the other hand, when the degrees of the denominator and of the numerator are the same, writing the rational interpolant in
a barycentric form allows to freely choose the weights appearing in this formula. These weights can be chosen by imposing various
additional conditions such as monotonicity or the absence of poles \cite{ber1,ber2,floa,ber3}.

For an interesting discussion between the coefficients of the interpolating rational function and the weights of its
barycentric representation, see \cite{berM}.
For the important problems of the ill--conditioning of rational interpolation, and of the numerical stability of the algorithms
for its solution, consult \cite{berM,peter}.

\vskip 2mm

In this paper, we will construct for the first time rational functions possessing both properties, that is interpolating $f$ at
some points of the complex plane, and whose series expansion around zero coincides with the Taylor series $f$ as far as possible.
Of course, this case is a particular instance of the general rational Hermite interpolation problem treated in its full generality
 in \cite{csww2}, for example.
Then, using a different number of conditions than required, we are able to construct rational interpolants in the
least squares sense.
We will also show how information on the poles and the zeros of $f$ could be included into these interpolants
in a style similar to the definition of the partial Pad\'{e} approximants \cite{ppa}.

\section{Problem statement}

We consider two different arguments.
\begin{itemize}
\item Let $f$ be a function whose Taylor series expansion around zero is known. A {\it Pad\'{e}--type
approximant} of $f$ is a rational function with an arbitrarily chosen denominator of
degree $k$, and whose numerator, also of degree $k$, is determined such that the
power series expansion of the approximant around zero coincides with the development of $f$ as far
as possible, that is up to the term of degree $k$ inclusively \cite{pta}.
By choosing the denominator appropriately, this rational function has a
series expansion which agrees with that of $f$ up to the term of degree $2k$ inclusively. It is then
called a {\it Pad\'{e} approximant}, and there is no freedom in the choice of the coefficients of the numerator and
the denominator of the rational approximant. On this topic, see, for example \cite{bak,birk}.
\item Let $f$ be a function whose values at $k+1$ distinct points in the complex plane are known. It is possible to
construct a rational function, with a numerator and a denominator both of degree $k$, which interpolates $f$
at these points. If this rational function is written in barycentric form, it depends on $k$ nonzero weights which can
be arbitrarily chosen. But, by Thiele's interpolation formula, it is also possible to obtain a rational function, with a numerator and
a denominator both of degree $k$ which interpolates $f$ at $2k+1$ distinct points in the complex plane. In that case, there is no freedom in the construction of the rational interpolant.
\end{itemize}

We now consider these two themes together and work in both directions in a different way.
Each of these choices leads to a different rational function whose series expansion agrees with that of $f$
as far as possible, and which interpolates $f$ at distinct points in the complex plane.

\begin{itemize}
\item We determine the denominator of the Pad\'{e}--type approximant so that it also interpolates $f$ at as many
distinct points in the complex plane as possible, that is $k$ points. Thus we obtain a rational function
interpolating $f$ at $k$ points and with an order of approximation $k+1$ at 0.
Such a rational function will be called a {\it Pad\'{e}--type rational interpolant}.
\item We determine the weights of the barycentric formula for the rational interpolant so that its power series
expansion coincides
with that of $f$ as far as possible, that is up to the term of degree $k-1$ inclusively.
This approach produces a rational function
with an order of approximation $k$ at 0, and interpolating $f$ at $k+1$ points.
Such a rational function will be called a {\it Pad\'{e}--type barycentric interpolant}.
\end{itemize}

In each case, different interpolation or approximation conditions can be considered, and the rational
function can be computed in the least squares sense.
Rational interpolants with arbitrary degrees in the numerator and in the denominator of the interpolant could also be defined similarly.
Let us mention that it is also possible to work with the reciprocal function $g$ of $f$, and its
reciprocal series which is defined by the algebraic relation $f(t)g(t)=1$.

\vskip 2mm

In the sequel, the formal power series $f$ will be written as
$$f(t)=c_0+c_1t+c_2t^2+\cdots$$

\section{Pad\'{e}--type rational interpolants}

We will begin by treating the case of a formal power series and, then, we will consider a series
in Chebyshev polynomials.

\subsection{Power series}

Let $R_k$ be written as
$$R_k(t)=\frac{N_k(t)}{D_k(t)}=\frac{a_0+a_1t+\cdots+a_kt^k}{b_0+b_1t+\cdots+b_kt^k}.$$
If the coefficients $b_i$ of the denominator are arbitrarily chosen (with $b_k \neq 0$), and if the coefficients $a_i$ of
the numerator are computed by the relations
\begin{equation}
\label{cp}
\left.
\begin{array}{rcl}
a_0&=&c_0b_0 \\
a_1&=&c_1b_0+c_0b_1\\
&\vdots&\\
a_k&=&c_kb_0+c_{k-1}b_1+\cdots+c_0b_k
\end{array}
\right\}
\end{equation}
then $R_k$ is the {\it Pad\'{e}--type approximant} $(k/k)_f$ of $f$ which satisfies the approx\-im\-at\-ion--through--order
conditions $f(t)-R_k(t)={\cal O}(t^{k+1})$. Let us remind that this condition means that $f$, $R_k$, and
their derivatives up to the $k$th inclusively take the same values at the point $t=0$.
Replacing $a_0,\ldots,a_k$ by their expressions (\ref{cp}) in $N_k$, and gathering the terms corresponding to each
$b_i$, we also have
\begin{equation}
N_k(t)=b_0S_k(t)+b_1tS_{k-1}(t)+\cdots+b_kt^kS_0(t), \label{numer}
\end{equation}
with
\begin{equation}
\label{psum}
S_n(t)=c_0+c_1t+\cdots+c_nt^n, \qquad n=0,1,\ldots
\end{equation}

Let us now determine $b_0,\ldots,b_k$ such that $R(\tau_i)=f(\tau_i)(=: f_i)$ for $i=1,\ldots,l$, that is such that
$$N_k(\tau_i)-f_iD_k(\tau_i)=0, \qquad i=1,\ldots,l,$$
where $\tau_1,\ldots,\tau_l$ are distinct points in the complex plane (none of them being 0).
We obtain the system
\begin{equation}
\label{sys}
(S_k(\tau_i)-f_i)b_0+\tau_i(S_{k-1}(\tau_i)-f_i)b_1+\cdots+\tau_i^k(S_0(\tau_i)-f_i)b_k=0, \quad i=1,\ldots,l.
\end{equation}
Since a rational function is defined up to a multiplying factor, we set $b_0=1$ (imposing another normalization
condition could lead to $b_0=0$ and, so, $a_0=0$, thus reducing the degree), and
we obtain a system of $l$ linear equations in the $k$ unknowns $b_1,\ldots,b_{k}$.
We consider its least squares solution if $l>k$ (overdetermined system), and its minimum norm solution for $l\leq k$ (underdetermined or singular system). The system has always a unique solution which determines a unique rational interpolant.
Therefore, the $b_i$'s are first determined by the interpolation conditions and, then, the $a_i$'s are calculated by
formulae (\ref{cp}).

\vskip 2mm

Multiplying each equation in (\ref{sys}) by $\tau_i^{-k}$ (the reason will be made clear later) and using (\ref{numer}),
we obtain the following Property, assuming that the denominator is different from zero.

\begin{property} \label{prop1}
~~\\
When $l=k$, it holds
$$R_k(t)=\frac{\left|
\begin{array}{cccc}
S_k(t) & tS_{k-1}(t) & \cdots & t^kS_0(t) \\
\tau_1^{-k}(S_k(\tau_1)-f_1) & \tau_1^{-k+1}(S_{k-1}(\tau_1)-f_1) & \cdots & S_0(\tau_1)-f_1 \\
\vdots & \vdots && \vdots \\
\tau_k^{-k}(S_k(\tau_k)-f_k) & \tau_k^{-k+1}(S_{k-1}(\tau_k)-f_k) & \cdots & S_0(\tau_k)-f_k
\end{array}
\right|}
{\left|
\begin{array}{cccc}
1 & t & \cdots & t^k \\
\tau_1^{-k}(S_k(\tau_1)-f_1) & \tau_1^{-k+1}(S_{k-1}(\tau_1)-f_1) & \cdots & S_0(\tau_1)-f_1 \\
\vdots & \vdots && \vdots \\
\tau_k^{-k}(S_k(\tau_k)-f_k) & \tau_k^{-k+1}(S_{k-1}(\tau_k)-f_k) & \cdots & S_0(\tau_k)-f_k
\end{array}
\right|}.
$$
\end{property}

\noindent{\it Proof:}~~\\
Let us take $t=\tau_i$ in this formula, and multiply the first row of the numerator and of the denominator by
$\tau_i^{-k}$. Then, subtract the row $i+1$ of the numerator from the first one. This first row becomes
$\tau_i^{-k}f_i, \tau_i^{-k+1}f_i,\ldots,f_i$, and we obtain $R_k(\tau_i)=f_i$, for $i=1,\ldots,k$, since the first row of
the denominator is $\tau_i^{-k}, \tau_i^{-k+1},\ldots,1$. Thus the interpolation property of $R_k$ has been recovered
from its determinantal expression. $\Box$

\vspace{0.2cm}

Let us now define the linear functionals $L_i$ acting on the vector space of polynomials by (this is the reason for multiplying
each equation in (\ref{sys}) by $\tau_i^{-k}$)
$$L_i(t^j)=\tau_i^{-j}(S_j(\tau_i)-f_i), \quad j=0,1,\ldots, \quad i=1,2,\ldots$$
The polynomial
$$P_k(t)=D_k
\left|
\begin{array}{cccc}
t^k & t^{k-1} & \cdots & 1 \\
\tau_1^{-k}(S_k(\tau_1)-f_1) & \tau_1^{-k+1}(S_{k-1}(\tau_1)-f_1) & \cdots & S_0(\tau_1)-f_1 \\
\vdots & \vdots && \vdots \\
\tau_k^{-k}(S_k(\tau_k)-f_k) & \tau_k^{-k+1}(S_{k-1}(\tau_k)-f_k) & \cdots & S_0(\tau_k)-f_k
\end{array}
\right|,$$
where $D_k$ is any nonzero normalization factor, satisfies the so--called {\it biorthogonality conditions}
$$L_i(P_k(t))=0, \quad i=1,\ldots,k, \quad L_{k+1}(P_k) \neq 0.$$
Such a polynomial is the $k$th member of the family of {\it formal biorthogonal polynomials} with respect to the linear
functionals $\{L_i\}$ \cite[pp. 104ff.]{biort}, and we see that the denominator of $R_k$ is equal to
$\widetilde P_k(t)=t^k P_k(t^{-1})$. This polynomial may not exist for some values of $k$, or its degree may be less than $k$. There is no general theory about that but, when it exists, $P_k$ is unique up to its normalization factor.

Let now $c$ be the linear functional acting on the vector space of polynomials and defined by
$c(x^i)=c_i$ for $i=0,1,\ldots$, let $Q_k$ be the polynomial of degree $k-1$ in $t$
$$Q_k(t)=c\left(x \frac{P_k(x)-P_k(t)}{x-t}\right),$$
and set $\widetilde Q_k(t)=t^{k-1}Q_k(t^{-1})$. From the definitions of $\widetilde P_k$, $\widetilde Q_k$,
and the determinantal formula of $R_k$ given in Property \ref{prop1}, we have the following Property.

\begin{property}
$$R_k(t)=c_0+t\frac{\widetilde Q_k(t)}{\widetilde P_k(t)},  \mbox{~~~~when~} l=k.$$
\end{property}

This Property shows that $R_k$ is exactly the generalization of the Pad\'{e}--type approximants defined in \cite[pp. 97ff.]{biort},
and, thus, it holds $R_k(t)-f(t)={\cal O}(t^{k+1})$ as required by our approximation--through--order conditions.

\vskip 2mm

It is possible to construct Pad\'{e}--type rational interpolants $(p/q)_f$ with an arbitrary degree $p$ in the numerator and
$q$ in the denominator, and then to determine its denominator in order to satisfy $q$ (or even $l \neq q$)
interpolation conditions \cite{hand,iser}.
Let us set $N_p(t) = a_0 + a_1t+ \cdots + a_pt^p$, and $D_q(t) = b_0 + b_1t + \cdots+ b_qt^q$.
The coefficients of the denominator are first computed as the solution of the system (\ref{sys})
with $l=q$ (or even $l \neq q$).
Then, the coefficients of the numerator are given by
\begin{equation*}
\left.
\begin{array}{rcl}
a_0 & = & c_0b_0 \\
a_1 & = & c_1b_0 + c_0b_1 \\
&\vdots &  \\
a_p & = & c_pb_0 + c_{p-1}b_1 + \cdots + c_{p-q}b_q,
\end{array}
\right\}
\end{equation*}
with the convention that $c_i = 0$ for $i < 0$, and the partial sums (\ref{psum}) computed accordingly. Such
an interpolant satisfies $(p/q)_f(\tau_i)=f_i$ for $i=1,\ldots,q$ and $(p/q)_f(t)-f(t)={\cal O}(t^{p+1})$.

\vskip 2mm

If some poles and some zeros of $f$ are known, this information could be included into the construction of the
rational interpolant. Let $p_1,\ldots,p_m$ and $z_1,\ldots,z_n$ be these poles and zeros, respectively.

Setting $P_m(t)=(t-p_1)\cdots(t-p_m)$ and $Z_n(t)=(t-z_1)\cdots(t-z_n)$, we are looking for the rational function
$$R_k(t)=\frac{N_k(t)Z_n(t)}{D_k(t)P_m(t)}$$
such that $R_k(\tau_i)=f(\tau_i)(=f_i)$ for $i=1,\ldots,k$, and such that
$f(t)-R_k(t)={\cal O}(t^{k+1})$. Such a rational function is called a {\it partial Pad\'{e}--type rational interpolant}
since it is similar to the partial Pad\'{e} approximants introduced
in \cite{ppa}, but with a lower order of approximation.

We must have
\begin{eqnarray*}
&&N_k(\tau_i)Z_n(\tau_i)-f_iD_k(\tau_i)P_m(\tau_i)=0 \\
&&N_k(\tau_i)-f_i\frac{P_m(\tau_i)}{Z_n(\tau_i)}D_k(\tau_i)=0, \qquad i=1,\ldots,k.
\end{eqnarray*}
Setting $N_k$ and $D_k$ as above, the coefficients of $D_k$ are first determined as the preceding ones with
$f_i$ replaced by $f_iP_m(\tau_i)/Z_n(\tau_i)$ in the system (\ref{sys}), and then the coefficients of $N_k$
are obtained by the same relations
as before where, now, the coefficients $c_i$ have to be replaced by those of the series expansion of $f(t)P_m(t)/Z_n(t)$
in (\ref{psum}).  Thus,
we first compute the coefficients of $h(t)=f(t)/Z_n(t)$ by identification in the relation $f(t)=h(t)Z_n(t)$. Then
the coefficients of $f(t)P_m(t)/Z_n(t)=h(t)P_m(t)$ are obtained by a simple product. These coefficients replace
the $c_i$'s in the definition of the partial sums (\ref{psum}).
Let us mention that this division and the following multiplication can be performed monomial by monomial
in order to avoid the computation of the coefficients of the polynomials $Z_n$ and $P_m$.  Indeed, we can begin
by computing the coefficients of $f(t)/(t-z_1)$, then, from these coefficients, we compute those
of $(f(t)/(t-z_1))/(t-z_2)$,
and so on until the division by $(t-z_n)$. Thus, we obtain the coefficients of $h$. Then, we formally multiply
$h(t)$ by $(t-p_1)$, the result by $(t-p_2)$, and so on until $(t-p_m)$ which gives the coefficients of
$h(t)P_m(t)=f(t)P_m(t)/Z_n(t)$.

\subsection{Fourier and Chebyshev series}

Fourier series can be approximated similarly by a procedure introduced in
\cite{pw} and developed in \cite{gibbs}. It consists in adding to the Fourier series its conjugate series, thus
transforming it, by a change of variable, into a power series, then computing the interpolants as
described above, and finally keeping only their real part. The approximation of parametric curves is another
topic which could be explored.

\vskip 2mm

Let us consider the case of a series in Chebyshev polynomials
$$f(t)=\frac{c_0}{2}+\sum_{i=1}^\infty c_i T_i(t),$$
where $T_i(t)=\cos( i \arccos t)$. The rational interpolant $R_k$ is defined as
$$R_k(t)=\frac{h_0/2+h_1T_1(t)+\cdots+h_kT_k(t)}{e_0/2+e_1T_1(t)+\cdots+e_kT_k(t)}.$$
Adapting to our case a general approach due to Hornecker \cite{horna,hornb} and particularized by
Paszkowski \cite{pasz} using the multiplication law
$T_i(t)T_j(t)=(T_{|i-j|}(t)+T_{i+j}(t))/2$ for Chebyshev polynomials, we have
$R_k(t)-f(t)={\cal O}(T_{k+1}(t))$ for any choice of the coefficients $e_i$ of the denominator,  if the coefficients
$h_i$ of the numerator are computed by
\begin{eqnarray*}
h_0 &=& c_0e_0/2+\sum_{i=1}^k c_ie_i\\
h_n &=& (c_ne_0+\sum_{j=1}^k (c_{|n-j|}+c_{n+j})e_j)/2, \quad n=1,\ldots,k.
\end{eqnarray*}
Let us now choose $e_0,\ldots,e_k$ such that $R_k(\tau_i)=f_i$ for $i=1,\ldots,k$. Similarly to the procedure followed
for a power series, these coefficients must satisfy
\begin{eqnarray*}
&&c_0e_0/2+\sum_{j=1}^k c_je_j+\sum_{n=1}^k \left(c_ne_0+\sum_{j=1}^k (c_{|n-j|}+c_{n+j})e_j\right)
T_n(\tau_i)-\\
&&\mbox{~~~~~~~~~~~~~~~~~~~~~~~~~~~~~~~~~~~~~~~~~~~~~~~~~~~~~~~~~~} f_i\left(e_0+2\sum_{j=1}^k e_jT_j(\tau_i)\right)=0,
\end{eqnarray*}
for $i=1,\ldots,k$, thus leading to the system
\begin{eqnarray*}
&&\left(c_0/2+\sum_{n=1}^kc_nT_n(\tau_i)-
f_i\right)e_0+\\
&& \mbox{~~~~~~~~~~~~~~~~~~~~~~~}\sum_{j=1}^k \left(c_j+
\sum_{n=1}^k (c_{|n-j|}+c_{n+j})T_n(\tau_i)-2f_iT_j(\tau_i)\right)e_j=0,
\end{eqnarray*}
for $i=1,\ldots,k$. Since a rational function is defined apart a multiplying factor, we set $e_0=1$ for solving it.

This approach can be extended to a numerator of degree $n+k$, $k \geq 1$ \cite[pp. 161ff.]{cont}, \cite[pp. 220ff.]{birk}. Moreover, since a Chebyshev series is a cosine series, its conjugate series could be added to it, as indicated above for Fourier series, and then a rational Pad\'{e}--type interpolant could be constructed, keeping only its real part.

\section{Pad\'{e}--type barycentric  interpolants}

We consider the following barycentric rational function 
$$R_k(t)=\frac{\displaystyle\sum_{i=0}^k\frac{w_i}{t-\tau_i}f_i}{\displaystyle\sum_{i=0}^k\frac{w_i}{t-\tau_i}},$$
where $f_i=f(\tau_i)$. This rational function interpolates $f$ at the $k+1$ points $\tau_i$, $i=0,\ldots,k$, whatever
the $w_i \neq 0$ are.
It is well--known that, by the Lagrangian interpolation formula for the denominator of $R_k$,
$w_i=q_i/v'(\tau_i)$ with $v(t)=\prod_{j=0}^k (t-\tau_j)$, and
$v'(\tau_i)=\prod_{j=0, j \neq i}^k (\tau_i-\tau_j)$,  where $q_i$ is the value of the denominator of $R_k$ at
the point $\tau_i$. This remark shows that, as in the case of Pad\'{e}--type rational interpolation,
the rational interpolant $R_k$ is fully determined by its denominator as mentioned in \cite{berM}.
Let us remind that, for the choice $w_i=v'(\tau_i)$, $R_k$ becomes a polynomial and that, for several choices
of the points $\tau_i$ closed expressions of the weights $w_i$ are known.

\vskip 2mm

Let us now determine $w_0,\ldots,w_k$ such that
$$f(t)-R_k(t)={\cal O}(t^{k}).$$
In that case, $R_k$ is a Pad\'{e}--type approximant $(k/k)_f$ of $f$, but with a lower order $k$ of approximation
instead of $k+1$.
This condition means that $f$ and $R_k$ and
their derivatives up to the $(k-1)$th inclusively take the same values at the point $t=0$.
Let us mention that it is not possible to improve the order of approximation for obtaining an exact
Pad\'{e}--type approximant.

The preceding approximation--through--order condition reads
$$\sum_{i=0}^k \frac{w_i}{t-\tau_i}f_i=(c_0+c_1t+\cdots)\sum_{i=0}^k \frac{w_i}{t-\tau_i}+{\cal O}(t^{k}).$$
Dividing each fractional term by the corresponding $\tau_i$ (obviously all the $\tau_i$ have to be different from zero,
which is not a restriction since our Pad\'{e}--type barycentric interpolant will interpolate $f$ at $t=0$), changing the signs,
and using the formal identity
$$\frac{1}{1-t/\tau_i}=1+\frac{t}{\tau_i}+\frac{t^2}{\tau_i^2}+\cdots,$$
we have
\begin{eqnarray*}
&&\sum_{i=0}^k\frac{w_i}{\tau_i}f_i\left(1+\frac{t}{\tau_i}+\frac{t^2}{\tau_i^2}+\cdots\right)=\\
&& \mbox{~~~~~~~~~~~~~~~~~~~~~~~~}(c_0+c_1t+\cdots)
\sum_{i=0}^k\frac{w_i}{\tau_i}\left(1+\frac{t}{\tau_i}+\frac{t^2}{\tau_i^2}+\cdots\right)+{\cal O}(t^{k}).
\end{eqnarray*}
Identifying the coefficients of identical powers of $t$ on both sides leads to
\begin{eqnarray*}
&&\sum_{i=0}^k\frac{w_i}{\tau_i}f_i=c_0\sum_{i=0}^k\frac{w_i}{\tau_i}\\
&&\sum_{i=0}^k\frac{w_i}{\tau_i}f_i\frac{1}{\tau_i}=c_0\sum_{i=0}^k\frac{w_i}{\tau_i}\frac{1}{\tau_i}+c_1\sum_{i=0}^k\frac{w_i}{\tau_i} \\
&&\sum_{i=0}^k\frac{w_i}{\tau_i}f_i\frac{1}{\tau_i^2}=c_0\sum_{i=0}^k\frac{w_i}{\tau_i}\frac{1}{\tau_i^2}+c_1\sum_{i=0}^k\frac{w_i}{\tau_i}
\frac{1}{\tau_i}+c_2\sum_{i=0}^k\frac{w_i}{\tau_i},
\end{eqnarray*}
and so on up to the term of degree $k-1$ inclusively.

Thus, the $w_i$ must be the solution of the linear system
\begin{equation}
\left.
\label{sysbar}
\begin{array}{lll}
&&\displaystyle \sum_{i=0}^k(f_i-c_0)\frac{w_i}{\tau_i}=0\\
&&\displaystyle \sum_{i=0}^k\left(\frac{f_i}{\tau_i}-\frac{c_0}{\tau_i}-c_1\right)\frac{w_i}{\tau_i}=0\\
&&\cdots \cdots \cdots \cdots \cdots \cdots \cdots \cdots \cdots \cdots \\
&&\displaystyle \sum_{i=0}^k\left(\frac{f_i}{\tau_i^{k-1}}-\frac{c_0}{\tau_i^{k-1}}-\frac{c_1}{\tau_i^{k-2}}-\cdots-c_{k-1}\right)\frac{w_i}{\tau_i}=0.
\end{array}
\right\}
\end{equation}
Since a rational fraction is defined apart a multiplying factor in its numerator and in its denominator, we will set $w_0=1$
and, thus, we obtain a system of $k$ equations in the $k$ unknowns $w_1,\ldots,w_k$.

This approach needs the knowledge of the values of $f$ at $k+1$ points, and that of the
coefficients $c_0,\ldots,c_{k-1}$.

\vskip 2mm

Let us write the system (\ref{sysbar}) as
$$\sum_{i=0}^k a_{ji}w_i=0, \quad j=1,\ldots,k.$$
Then, we obtain two determinantal expressions for $R_k$, the first one in a barycentric form, and the second
one in a Lagrangian--type basis.

\begin{property}
$$R_k(t)=\frac{
\left|
\begin{array}{ccc}
f_0/(t-\tau_0) & \cdots & f_k/(t-\tau_k)\\
a_{10} & \cdots & a_{1k} \\
\vdots && \vdots \\
a_{k0} & \cdots & a_{kk}
\end{array}
\right|}
{\left|
\begin{array}{ccc}
1/(t-\tau_0) & \cdots & 1/(t-\tau_k)\\
a_{10} & \cdots & a_{1k} \\
\vdots && \vdots \\
a_{k0} & \cdots & a_{kk}
\end{array}
\right|}
=
\frac{
\left|
\begin{array}{ccc}
f_0L_0(t) & \cdots & f_kL_k(t)\\
a_{10} & \cdots & a_{1k} \\
\vdots && \vdots \\
a_{k0} & \cdots & a_{kk}
\end{array}
\right|}
{
\left|
\begin{array}{ccc}
L_0(t) & \cdots & L_k(t)\\
a_{10} & \cdots & a_{1k} \\
\vdots && \vdots \\
a_{k0} & \cdots & a_{kk}
\end{array}
\right|},
$$
with, for $i=0,\ldots,k$,
\begin{eqnarray*}
a_{1i} &=& (f_i-c_0)\frac{w_i}{\tau_i}\\
a_{ji} &=& \frac{a_{j-1,i}-c_{j-1}{w_i}}{\tau_i}, \quad j=2,\ldots,k,
\end{eqnarray*}
and
$$L_i(t)=\prod_{\genfrac{}{}{0pt}{}{j=0}{j \neq i}}^k (t-\tau_j).$$
\end{property}

\noindent{\it Proof:}~~\\
The second formula comes out from $L_i(t)=L(t)/(t-\tau_i)$ with $L(t)=(t-\tau_0)\cdots(t-\tau_k)$.
Since $L_i(\tau_m)=0$ for $m \neq i$ and $L_i(\tau_i) \neq 0$, we immediately recover, from the second expression,
the interpolation property $R_k(\tau_i)=f_i$ for $i=0,\ldots,k$.
For recovering the approximation--through--order property, the expressions
$1/(t-\tau_i)$ in the numerator and in the denominator of $R_k$ have to be replaced by
$-1/(\tau_i(1-t/\tau_i))=-(1+t/\tau_i+t^2/\tau_i^2+\cdots)/\tau_i$, and the coefficient of each power of $t$
has to be separately identified up to the degree $k-1$ inclusively. $\Box$

\vspace{0.2cm}

Assume now that only $c_0,\ldots,c_{l-1}$ are known, with $l<k$.
We can choose $w_0,\ldots,w_k$ such that $f(t)-R_k(t)={\cal O}(t^l)$ by considering only the
first $l$ equations in the preceding system, and replacing the last ones by the equations
$$\sum_{i=0}^k\left(\frac{f_i}{\tau_i^{l+j-1}}-\frac{c_0}{\tau_i^{l+j-1}}-\frac{c_1}{\tau_i^{l+j-2}}-\cdots-c_{l-1}\right)
\frac{w_i}{\tau_i}=0, \qquad j=1,\ldots,k-l,$$
which is equivalent to considering that the coefficients $c_l,\ldots,c_{k-1}$ are zero in the system (\ref{sysbar}).
The rational function $R_k$ now interpolates $f$ in $k+1$ points and its expansion coincides with that of $f$
up to the term of degree $l-1$ inclusively.

\vskip 2mm

It is also possible to consider the case where $l>k$ coefficients of the series of $f$ are known.
Adding to the preceding system the equations
$$\sum_{i=0}^k\left(\frac{f_i}{\tau_i^{j}}-\frac{c_0}{\tau_i^{j}}-\frac{c_1}{\tau_i^{j-1}}-\cdots-c_{j}\right)
\frac{w_i}{\tau_i}=0, \quad j=k,\cdots,l-1,$$
and solving it in the least squares sense leads to an approximation $R_k$ whose series expansion agrees with that
of $f$ only in the least squares sense, and which interpolates $f$ at $k+1$ points.

\vskip 2mm

Let us again consider the case where some poles and some zeros of $f$ are known. The rational function
$$R_k(t)\frac{Z_n(t)}{P_m(t)}=\frac{\displaystyle\sum_{i=0}^k\frac{w_i}{t-\tau_i}f_iP_m(\tau_i)/Z_n(\tau_i)}
{\displaystyle\sum_{i=0}^k\frac{w_i}{t-\tau_i}},$$
interpolates $f$ at the $k+1$ points $\tau_i$, $i=0,\ldots,k$, whatever the $w_i \neq 0$ are, and it has the poles
$p_1,\ldots,p_m$ and the zeros $z_1,\ldots,z_n$. Thus it can be constructed as above after replacing everywhere
$f_i$ by $f_iP_m(\tau_i)/Z_n(\tau_i)$, and we obtain $f(t)-R_k(t)Z_n(t)/P_m(t)={\cal O}(t^k)$.

When the poles of $f$ are known, an explicit expression for the weights of the near--best rational interpolants
in a Chebyshev sense can be obtained \cite{deun}. As mentioned in this paper, the knowledge of the poles
dramatically improves the interpolation process as can be seen from the numerical examples given there,
and also in \cite{ppa,deun2}.

\section{Study of the error}

Let us set $R_k(t)=N_k(t)/D_k(t)$ either for the Pad\'{e}--type rational interpolants or the Pad\'{e}--type barycentric interpolants. We have
$$f(t)D_k(t)-N_k(t)={\cal O}(t^n)$$
with $n=k+1$ in the first case and $n=k$ in the second one.

We assume that all the interpolations points $\tau_i$ are real and belong to an interval $[a,b]$ and that, in this interval, $f$ has poles $\alpha_1,\ldots,\alpha_\nu$ of respective multiplicities $r_1,\ldots,r_\nu$ with $r_1+\cdots+r_\nu=m \leq n-1$, that none of these poles coincides with an interpolation point, and that, outside of the poles, $f$ has a bounded $(n+k)$th derivative. We set
$$\phi(t)=(t-\alpha_1)^{r_1}\cdots(t-\alpha_j)^{r_\nu}.$$
Thus, $f(t)\phi(t)$ is bounded in $[a,b]$. Let $\psi$ be a polynomial such that $Q(t)=\phi(t)\psi(t)$ has degree $n-1$. We write the error under the form
$$f(t)-R_k(t)=g(t)\frac{t^n(t-\tau_1)\cdots(t-\tau_k)}{D_k(t)Q(t)},$$
and set
$$w(x)=f(x)D_k(x)Q(x)-N_k(x)Q(x)-g(t)x^n(x-\tau_1)\cdots(x-\tau_k), \quad t \in [a,b].$$
The function $w$ has a simple zero at $x=t$ (by definition of the error), a zero of multiplicity $n$ at $x=0$, and simple zeros at $x=\tau_1,\ldots,\tau_k$.
Therefore, by Rolle's theorem and since the $(n+k)$th derivative of $N_k(x)Q(x)$ is identically zero, there exists a point $\xi_t \in [a,b]$, which depends on $t$, such that $w^{(n+k)}(\xi_t)=0$. Thus
$$g(t)=\frac{1}{(n+k)!}\frac{d^{n+k}}{d\xi^{n+k}}[f(\xi)D_k(\xi)Q(\xi)]\Big|_{\xi=\xi_t},$$
and it follows

\begin{property}
~~\\
Under the preceding assumptions
$$f(t)-R_k(t)=\frac{t^n(t-\tau_1)\cdots(t-\tau_k)}{(n+k)!D_k(t)Q(t)}
\frac{d^{n+k}}{d\xi^{n+k}}[f(\xi)D_k(\xi)Q(\xi)]\Big|_{\xi=\xi_t}, \quad t, \xi_t \in [a,b].$$
\end{property}

If $f$ has no pole in $[a,b]$, one can take $Q(x)=D_k(x)$ when $n=k+1$.
This result is adapted from \cite[pp. 116--7]{milne}.

\section{Numerical examples}

We will now show some numerical examples which gather several interesting
properties that will allow us to  exemplify the effectiveness of our procedures. But, before, let us give the
following {\it consistency property}

\begin{property}
~~\\
If $f$ is a rational function with a numerator and a denominator both of degree smaller or equal to $k$, then,
our two procedures produce a rational function $R_k$ which is identical to $f$ when $l=k$.
\end{property}

\noindent{\it Proof:}
~~\\
This property comes out from the fact that $R_k$ is defined by a set of linear equations which  is the same
as the set of equations which defines $f$, and the result follows from the uniqueness of $R_k$. $\Box$

\vspace{0.2cm}

Our numerical experiments were performed using Matlab$^\circledR$ 7.11.
Let us remind that the solution of a rectangular system of equations $Ax=b$ of maximal rank $r=\min(l,k)$
with $A \in {\mathbb{C}}^{l \times k}$
is $x=A^\dag b$, where $A^\dag$ is the Moore--Penrose pseudo--inverse of $A$ defined by
$A^\dag=(A^*A)^{-1}A^*$ if $r=k \leq l$
(overdetermined system) and $A^\dag=A^*(AA^*)^{-1}$ if $r=l \leq k$ (undertermined or singular system).
If the rank $r$ is not maximal, then $A^\dag=V\Sigma^\dag U^*$ where $A=U\Sigma V^*$ is the singular value
decomposition of $A$. The Matlab$^\circledR$ instruction {\tt pinv(A)*b} gives the least squares solution when the system is
overdetermined (that is the unique solution minimizing the 2--norm of the residual if the matrix is full rank, and the vector of
minimal 2--norm among those minimizing the 2--norm of the residual, if not), and the minimal 2--norm solution
when the system is underdetermined or singular. In all cases, the computation is based on the singular value decomposition of $A$.

\vskip 2mm

All curves (except in Figure \ref{cosRAT}) represent the errors in logarithmic scale.

\begin{figure*}
 \includegraphics[width=0.75\textwidth]{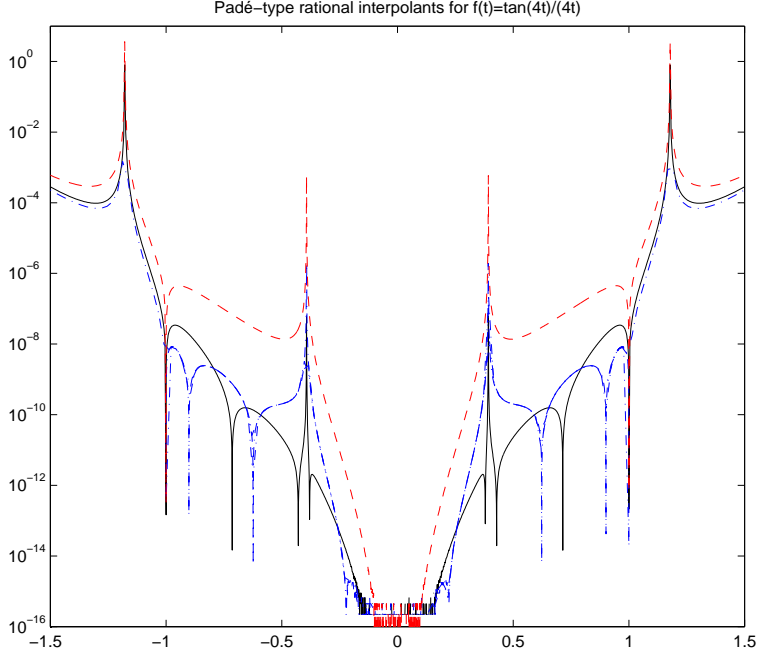}
\caption{Pad\'{e}--type rational interpolants with $k=8$ for $\tan(4t)/(4t)$: equidistant points in the interval $[-1,+1]$ (solid),
roots of unity (dashed), Chebyshev zeros (dash-dotted).}
\label{tanPT}
\end{figure*}

\begin{figure*}
\includegraphics[width=0.75\textwidth]{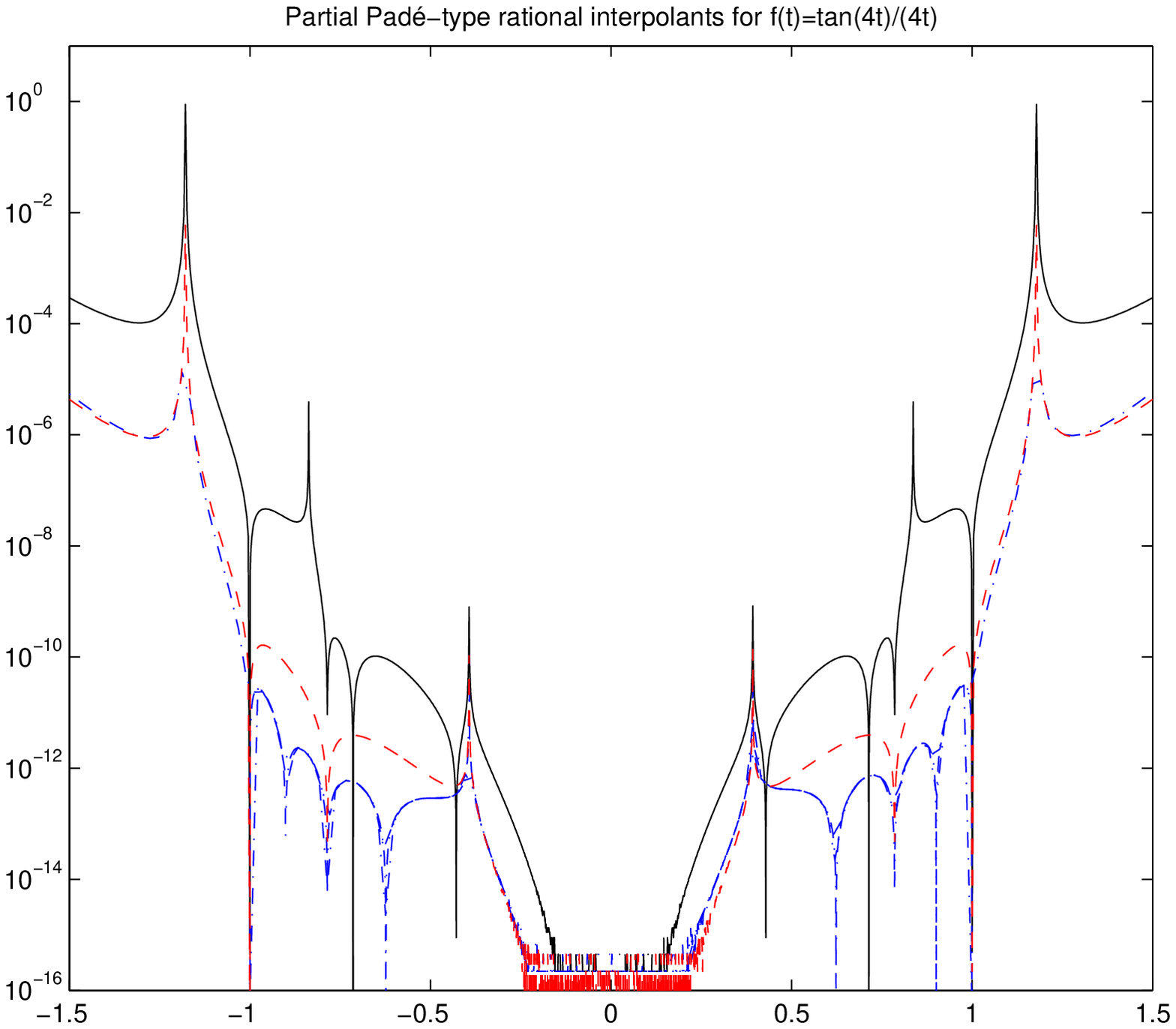}
\caption{Partial Pad\'{e}--type rational interpolants with $k=8$ for $\tan(4t)/(4t)$: equidistant points in the interval $[-1,+1]$ (solid),
roots of unity (dashed), Chebyshev zeros (dash-dotted).}
\label{tanPPT}
\end{figure*}

\subsubsection*{\bf Example 1: a function with poles}

 We consider the following function, and its series expansion
$$f(t)=\frac{\tan(\omega t)}{\omega t}=1+\frac{1}{3}\omega^2t^2+\frac{2}{15}\omega^4t^4
+\frac{17}{315}\omega^6t^6+\frac{62}{2835}\omega^8t^8+\cdots$$
This function has poles at odd multiples of $\pi/(2\omega)$, and zeros at odd multiples of $\pi/\omega$, except at 0.

\vskip 2mm

We considered three sets of interpolations points: equidistant points in the interval $[-1,+1]$,
the roots of unity, and the zeros of the Chebyshev polynomials of the first kind.
The complex choice was discussed in \cite{gonn}.
Let us mention that none of the interpolation points $\tau_i$ should be 0, since it is the point where
the Pad\'{e}--type approximants are computed and thus it always appears as an interpolation point.

\vskip 2mm

\noindent {\it Pad\'{e}--type rational interpolants}

The errors obtained with the Pad\'{e}--type rational interpolants are given in Figure \ref{tanPT} for $\omega=4$ and $k=8$.
The solid line corresponds to the real interpolation points, while the dashed one refers to the points
on the unit circle, and the dash-dotted one to the zeros of the Chebyshev polynomial.
The poles of $f$ are, in the interval considered in Figure \ref{tanPT}, at the points $\pm 0.39269908\ldots$ and $\pm 1.1780972\ldots$, and the zeros at $\pm 0.78539816\ldots$

Since the poles and the zeros are known, we took $Z_2(t)=(t-\pi/4)(t+\pi/4)$ and $P_2(t)=(t-\pi/8)(t+\pi/8)$.
The errors obtained with the partial Pad\'{e}--type rational interpolants are displayed in Figure \ref{tanPPT}.
Let us mention that, for some values of $k$, we could observe Froissart's doublets (nearby poles and zeros) that can be removed by the technique described below (Example 4).

The improvement brought by partially taking into account the knowledge of the poles and of the zeros is clear.
Choosing the zeros of the Chebyshev polynomials as the real interpolation points in $[-1,+1]$ does not change much
the quality of the results for such small values of $k$.

\vskip 2mm

\noindent {\it Pad\'{e}--type barycentric interpolants}

\begin{figure*}
\includegraphics[width=0.75\textwidth]{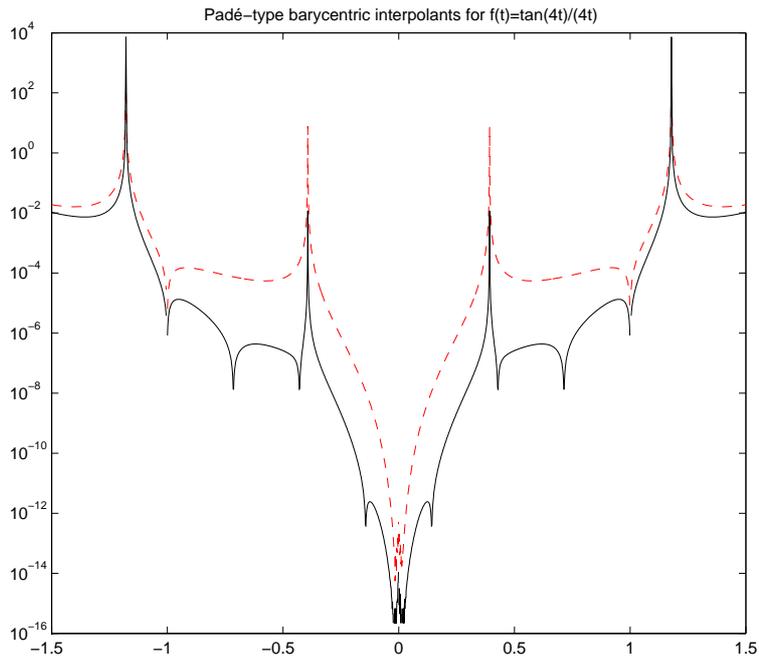}
\caption{Pad\'{e}--type barycentric interpolants with $k=8$ for $\tan(4t)/(4t)$: equidistant points in the interval $[-1,+1]$ (solid),
roots of unity (dashed).}
\label{tanBARY}
\end{figure*}

\begin{figure*}
\includegraphics[width=0.75\textwidth]{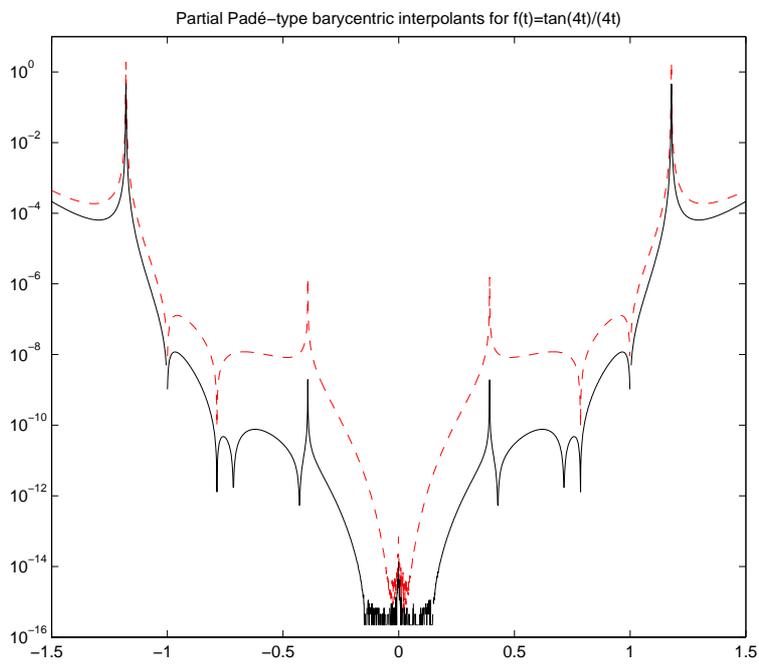}
\caption{Partial Pad\'{e}--type barycentric interpolants with $k=8$ for $\tan(4t)/(4t)$: equidistant points in the interval $[-1,+1]$ (solid),
roots of unity (dashed).}
\label{tanBARYp}
\end{figure*}

Let us now consider the same example but with the Pad\'{e}--type barycentric interpolants.
The results are given in Figure \ref{tanBARY}.
With the partial Pad\'{e}--type barycentric interpolants, we obtain the results of Figure \ref{tanBARYp}. In both figures, the solid line
corresponds to the real interpolation points, and the dashed one to the interpolation points on the unit circle.

Let us mention that, with the Shepard's weights $w_i=1/(t-\tau_i)$ \cite{shep}, the interpolants have poles around $-0.4$ and $+0.4$.

\subsubsection*{\bf Example 2: a function with a cut}

We consider the series
$$f(t)=\frac{\log(1+t)}{t}=1-\frac{t}{2}+\frac{t^2}{3}-\frac{t^3}{4}+\cdots$$
which converges in the unit disk and on the unit circle except at the point $-1$ since there is a cut from
$-1$ to $-\infty$.

\vskip 2mm

\noindent {\it Pad\'{e}--type rational interpolants}

\begin{figure*}
\includegraphics[width=0.75\textwidth]{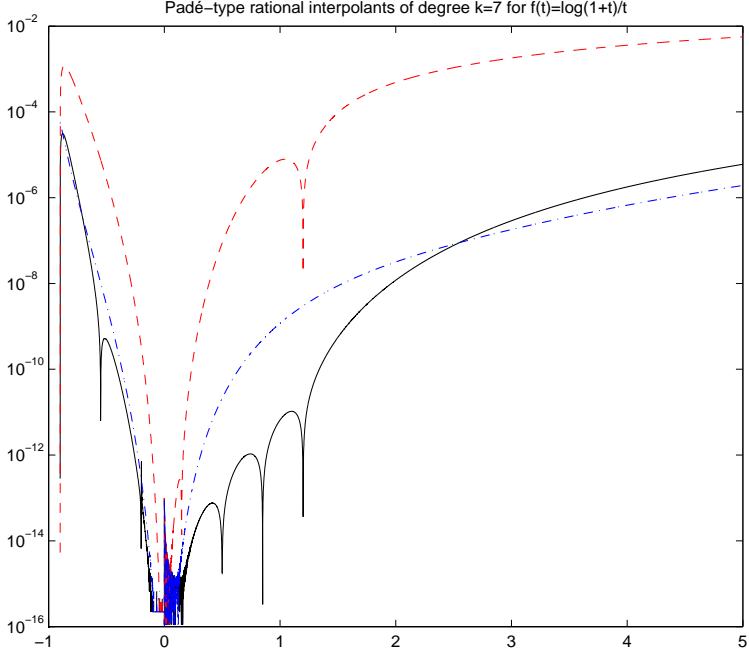}
\caption{Pad\'{e}--type rational interpolants with $k=7$ for $\log(1+t)/t$, and 3 (dashed), 7 (solid) and 14 (dash-dotted) interpolation points.}
\label{lslog}
\end{figure*}

For a Pad\'{e}--type interpolant of degree 7, we consider equidistant real interpolation points in the interval $[-0.9,+1.2]$.
For 7 points (solid line), the system to be solved is square. For 3 points (dashed line) and 14 points
(dash--dotted line), the system is solved in the least squares sense as explained above. The results are given in Figure \ref{lslog}.
We see that they are quite good even for values of $t$ far outside the convergence interval.

\vskip 2mm

\noindent {\it Pad\'{e}--type barycentric interpolants}

\begin{figure*}
\includegraphics[width=0.75\textwidth]{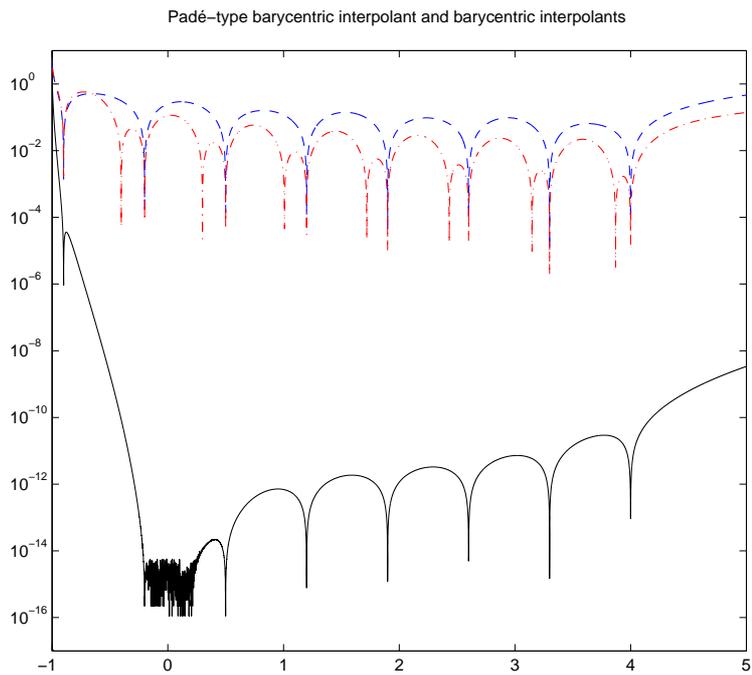}
\caption{Pad\'{e}--type barycentric interpolants with $k=7$ for $\log(1+t)/t$ (solid). Barycentric interpolants with Berrut weights (dashed), and Shepard weights (dash-dotted).}
\label{logBARY}
\end{figure*}

The interpolation points $\tau_i$ are taken equidistant in $[-0.9,+4]$, and $k=7$.
In Figure \ref{logBARY}, three types of weights $w_i$ are considered:
those corresponding to the Pad\'{e}--type barycentric interpolants are the same as explained above (solid line),
the weights $w_i=(-1)^i$ of Berrut \cite{ber2} (dashed line), and the weights $w_i=1/(t-\tau_i)$ suggested by Shepard
\cite{shep} (dash--dotted line), these last two choices ensuring pole--free interpolants on the real line.

\begin{figure*}
\includegraphics[width=0.75\textwidth]{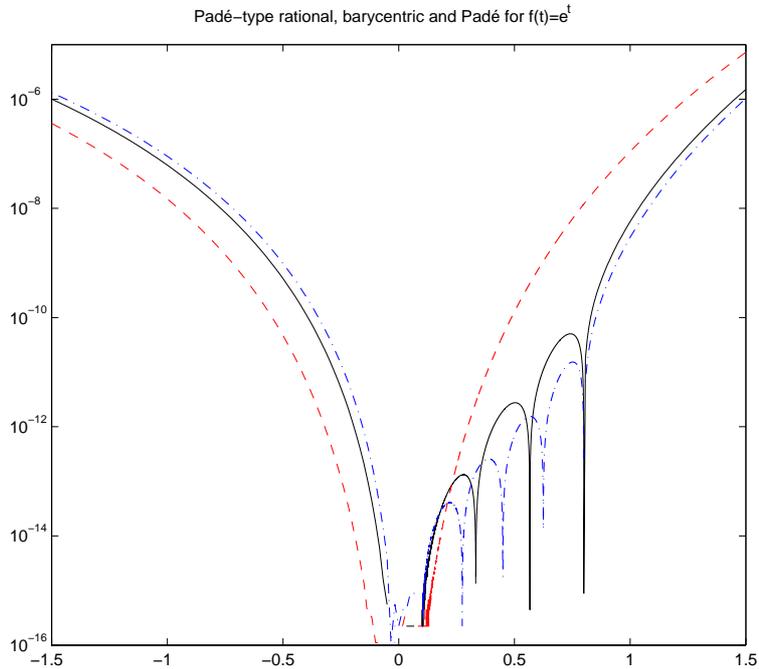}
\caption{Pad\'{e}--type rational (solid) and barycentric (dash-dotted) interpolants, and Pad\'{e} approximant (dashed) for $e^t$.}
\label{expRBP}
\end{figure*}

\subsubsection*{\bf Example 3: a continuous function}

We consider the exponential function
$$f(t)=e^t=1+\frac{t}{1!}+\frac{t^2}{2!}+\cdots$$
Let us now compare, for the degree $k=4$, the Pad\'{e}--type rational interpolant, the Pad\'{e}--type barycentric
interpolant, and the Pad\'{e} approximant $[4/4]$ which is given by
$$[4/4]_f(t)= (1680+840t+180t^2+20t^3+t^4)/(1680-840t+180t^2-20t^3+t^4).$$

Let us remind that $[4/4]_f(t)-e^t={\cal O}(t^9)$, and that its construction makes use of the first 8 coefficients
of the power series.
The results are given in Figure \ref{expRBP}, where the solid line represents the error of the
Pad\'{e}--type rational interpolant, the dashed line corresponds to the Pad\'{e} approximant, and the dash--dotted
line to the Pad\'{e}--type barycentric interpolant.

The interpolation points were chosen equidistant in the interval $[0.1,0.8]$. Notice that, for the interpolants,
the error is smaller around the interpolation points, while the errors of the Pad\'{e} approximant is more
symmetric around the origin.

\begin{figure*}
\includegraphics[width=0.75\textwidth]{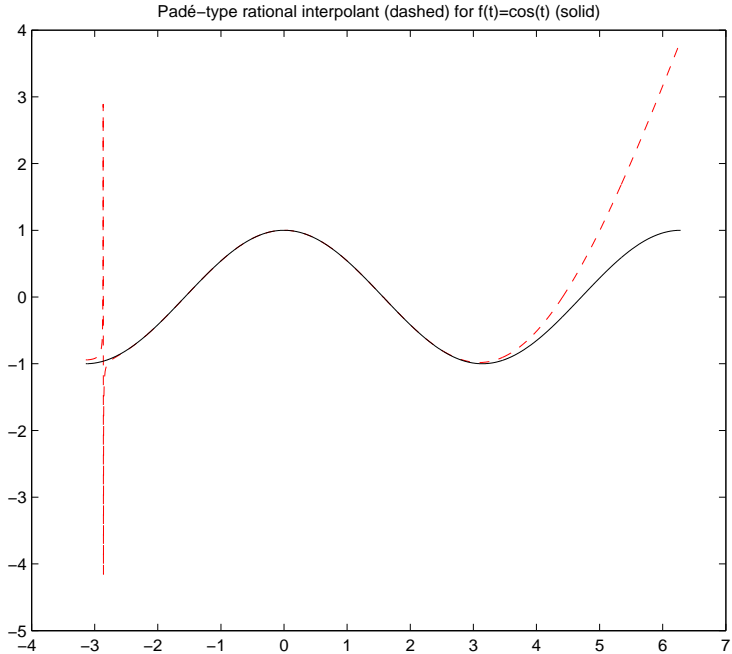}
\caption{Pad\'{e}--type rational interpolant (dashed) of the cosine function (solid).}
\label{cosRAT}
\end{figure*}

\subsubsection*{\bf Example 4: spurious pole removal}

Let us now give an example showing that the rational interpolant can have poles even if the function is continuous.
In fact, it is known \cite{csww1} that if, after cancelation of common factors between the numerator and the
denominator and ordering the interpolation points, two consecutive weights $w_i$ and $w_{i+1}$ in the barycentric formula
have the same sign, then the reduced interpolant has an odd number of poles in $[\tau_i,\tau_{i+1})$.

We consider the Pad\'{e}--type rational interpolant of the cosine function with 5 equidistant interpolation points in the
interval $[-\pi/2,+\pi/8]$. As may be seen in Figure \ref{cosRAT}, the interpolant (dashed line) has one
real pole at $t=-2.8636\ldots$ (its other poles are complex). When $t$ goes to infinity, the interpolant tends
to $25.269\ldots$

However, the results are quite good (the cosine function is the solid line in Figure \ref{cosRAT}) from the right
of the pole up to almost $\pi$.

\begin{figure*}
\includegraphics[width=0.75\textwidth]{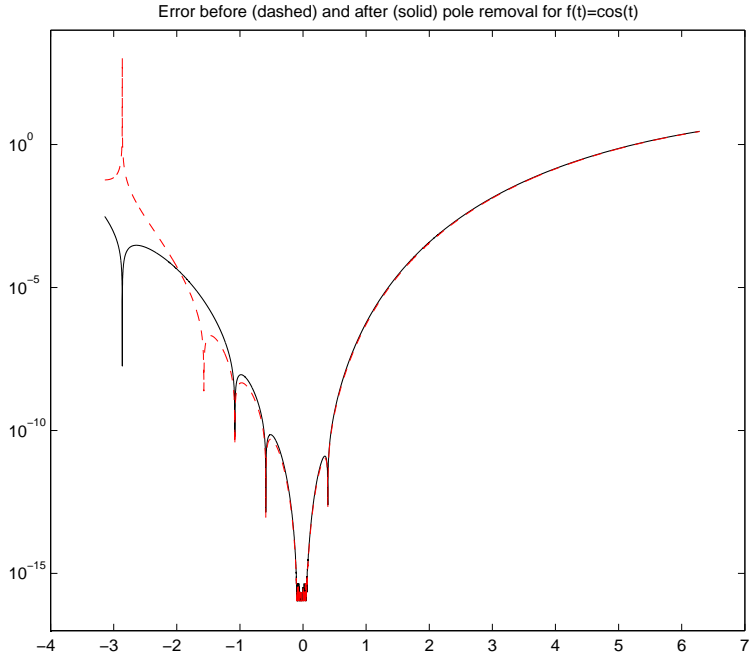}
\caption{Error before (solid) and after (dashed) the pole removal in the Pad\'{e}--type rational interpolant for $f(t)=\cos t$.}
\label{polerem}
\end{figure*}

It is possible to remove a spurious pole $p$ by forcing the Pad\'{e}--type interpolant to go through the point
$(p,f(p))$. In Figure \ref{polerem}, the first of the equidistant interpolation points is replaced by the pole $p$, a
procedure which removes it and leads to a better result (dashed line).

If the interpolant exhibits several poles, they can be eliminated successively.
If a new pole is introduced during the procedure, then it can be removed similarly.

In our case, the location of the pole was directly computed from the coefficients $b_i$ of the denominator
of the interpolant since they were available. It is also possible to locate approximately a pole when
the absolute value of the interpolant
becomes larger than a fixed threshold, or when the interpolant has a sudden change of sign,
and then to impose it as an interpolation point.

This procedure was tried on the Pad\'{e}--type barycentric interpolant for $\cos t$ in the same interval as
before, but with $k=13$. The interpolant was computed at 500 points in $[-\pi,+2\pi]$.
A sudden change of sign was observed in the interval $[5.3010,5.3199]$. We had
$R_{13}(5.3010)=15.068$ and $R_{13}(5.3199)=-25.101$. Replacing the first interpolation
point $\tau_0=-\pi/2$ by $\tau_0=5.3050$, the spurious pole was removed, and no other pole appeared.

The advantage of
this procedure is that it can also be used for Pad\'{e}--type barycentric interpolants where the coefficients
of the denominator are not explicitly known.

The same techniques can be applied to the case of partial Pad\'{e}--type rational and barycentric interpolation.

\section{Applications}

Let us now briefly discuss some possible applications to numerical analysis problems.

\subsection{Convergence acceleration}

We consider the sequence $(S_n=f(\tau_n))$ where $(\tau_n)$ is a sequence of parameters such that
$\lim_{n \to \infty}\tau_n=\tau_\infty \neq 0, \pm \infty$, and where $f$ is a function whose first coefficients
of the series expansion around 0 are known. We set $S=\lim_{n \to \infty} S_n$.

The convergence of the sequence $(S_n)$ can be accelerated by computing the Pad\'{e}--type rational interpolant
or the Pad\'{e}--type
barycentric interpolant $R_k^{(n)}$ satisfying $R_k^{(n)}(\tau_i)=S_i$ for $i=0,\ldots,k-1$ (or $k$ in the second
case), and setting $T_k^{(n)}=R_k^{(n)}(\tau_\infty)$.This is the essence of an extrapolation method
for accelerating the convergence of a sequence \cite{cbmrz}.
Under certain assumptions, the sequences $(T_k^{(n)})$ converge to $S$ faster than $(S_n)$ either when $k$ is
fixed and $n$ goes to infinity, or vice versa.

\subsection{Inversion of the Laplace transform}

We consider the Laplace transform
$$F(p)=\int_0^\infty e^{-p s}f(s) \, ds.$$

Assume that $F$ is known at some points $p_n$ for $n=0,\ldots,k-1$ (or $k$), and also the first coefficients of its series
expansion around 0. $F$ can be approximated by a Pad\'{e}--type rational interpolant or by a Pad\'{e}--type
barycentric interpolant $R_k$, and the interpolant then inverted, thus leading to an approximation of $f$.
Let us remark that, since $\lim_{p \to \infty} F(p)=0$, the degree of the numerator
of the interpolant must be smaller than the degree of its denominator.
The inversion can be performed without decomposing $F$ into its partial fractions by a procedure due to
Longman and Sharir \cite{long}.
Let $F$ have the form
$$F(p)=A \frac{p^{m}+\alpha_1 p^{m-1}+\cdots+\alpha_{m}}{p^n+\beta_1 p^{n-1}+\cdots+\beta_n},$$
with $m<n$. They showed that
$$f(s)=A \sum_{i=0}^\infty \frac{v_i}{i!}s^i$$
with
$$v_i = u_{i+m}+\alpha_1 u_{i+m-1}+\cdots+\alpha_{m}u_i, \quad i=0,1,\ldots,$$
where
\begin{eqnarray*}
u_i &=& 0, \qquad i=0,\ldots,n-2,\\
u_{n-1} &=& 1,\\
u_i &=& -(\beta_1 u_{i-1}+\cdots+\beta_n u_{i-n}), \quad i=n,n+1,\ldots
\end{eqnarray*}

Usually, the series giving $f$ is quickly converging.

\vskip 2mm

\begin{figure*}
\includegraphics[width=0.75\textwidth]{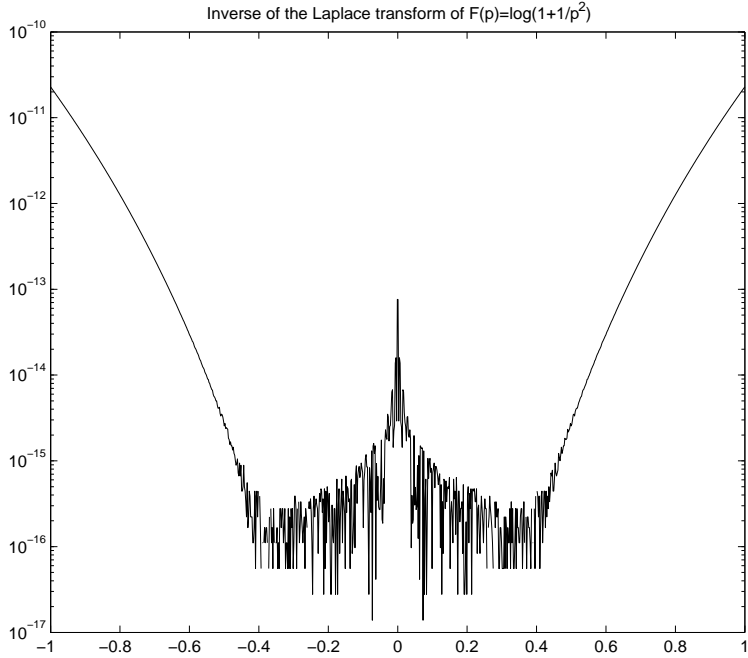}
\caption{Error for the inversion of the Laplace transform of $F(p)=\log(1+1/p^2)$.}
\label{lapl}
\end{figure*}

Let us take the example considered in \cite[p. 350]{cbmrz}
$$F(p)=\log(1+a^2/p^2), \qquad f(s)=2(1-\cos as)/s.$$
We make the change of variable $t=a^2/p^2$, and we set
$$F(p)=G(t)=\log(1+t)=t-\frac{t^2}{2}+\frac{t^3}{3}-\cdots$$

The Pad\'{e}--type (rational and barycentric) interpolants will be approximations of $G$.
Replacing $t$ by $a^2/p^2$ in a Pad\'{e}--type rational interpolant of degree $k$ in $t$ produces an
interpolant of degree $2k$ in $p$, and we obtain an approximation of $F$ of the form
$$R_{2k}(p)=\frac{a_0p^{2k}+a_1a^2p^{2k-2}+\cdots+a_ka^{2k}}{b_0p^{2k}+b_1a^2p^{2k-2}+\cdots+b_ka^{2k}}.$$
Notice that, since $c_0=0$ in the series expansion of $G(t)$, the relations (\ref{cp}) lead
to $a_0=0$, and, thus, $\lim_{p \to \infty} R_{2k}(p)=0$ which is consistent with the asymptotic property of the
Laplace transform. Thus, this approximant can be written as
$$R_{2k}(p)=A \frac{p^{2k-2}+\alpha_2 p^{2k-4}+\cdots+\alpha_{2(k-1)}}
{p^{2k}+\beta_2 p^{2k-2}+\cdots+\beta_{2k}},$$
with $A=a^2a_1/b_0$, $\alpha_{2i}=a^{2i}( a_{i+1}/a_1)$, for $i=1,\ldots,k-1$, $\beta_{2i}=a^{2i}(b_i/b_0)$, for
$i=1,\ldots,k$, the $\alpha$'s and the $\beta$'s with an odd index being zero.
We see that the series expansion of $R_{2k}(p)$ only
contains even powers of $1/p$ as the series $F(p)$ itself. Inverting $R_{2k}$ by the procedure of
Longman and Sharir \cite{long} (after replacing $m$ by $2k-2$ and $n$ by $2k$ in the formulae for
the $v_i$'s and the $u_i$'s),
or performing its partial fraction decomposition, gives an approximation of $f$.

\vskip 2mm

For $k=5$, $a=1$, $\tau_i=1/p_i^2$ with $p_i=0.1+ih$ for $i=0,\ldots,k-1$, and $h=2/(k-1)$,
the Pad\'{e}--type rational interpolant leads to the results of Figure \ref{lapl}, using 12 terms in the series
expansion of $f$. Although $f(0)=0$ and the series expansion by the method of Longman and Sharir is also 0
at $s=0$ (since $v_0=0$), there is a loss of accuracy around this point due to the indeterminacy.
These results have to be compared with those given in \cite[p. 350]{cbmrz} which were obtained
by constructing a rational interpolant with a numerator of degree 7 and a denominator of degree 8, that is using
16 interpolation points. We see that our Pad\'e--type rational interpolant
provides a much better precision. Moreover, the precision can be even improved by taking more terms in the series
for $f$ at almost no additional price.

This example could also be treated by making the change of variable $t=a/p$, thus leading to
$F(p)=G(t)=\log(1+t^2)=t^2-t^4/2+t^6/3-\cdots$.

\subsection{Piecewise rational interpolation}

Our approach can be used for constructing piecewise rational interpolants. Let $a<a' \leq 0 \leq b'<b$.
We construct a first Pad\'{e}--type rational or barycentric interpolant in $[a,a']$, and then a second one in $[b',b]$.
Due to the Pad\'{e}--type property of these interpolants and the fact that, for all $i$, $c_i=f^{(i)}(0)/i!$,
the two interpolants and their
first derivatives will have the same values at the point $t=0$. Obviously, by a change of variable, the same construction holds at a point different from the
origin, and it can be repeated.

One of the advantages of such a construction is to obtain a good accuracy with a low degree in the interpolants,
thus avoiding the usually bad conditioning when using more interpolation points and a rational interpolant with
a higher degree.

\begin{figure*}
\includegraphics[width=0.75\textwidth]{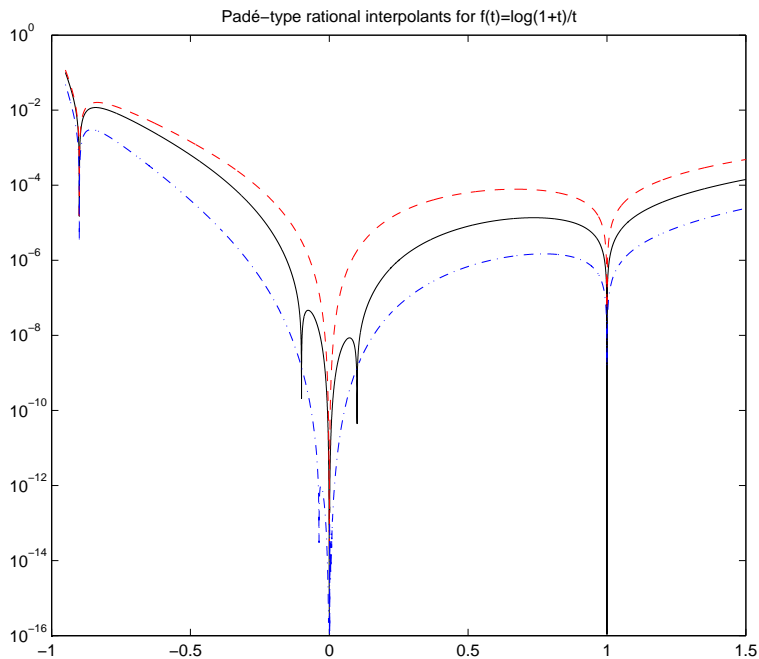}
\caption{Pad\'{e}--type rational interpolants for $\log(1+t)/t$: the piecewise case. }
\label{piecew}
\end{figure*}

\vskip 2mm

We interpolate the function $f(t)=\log(1+t)/t$ on the intervals $[-0.9,-0.1]$ and $[+0.1,+1]$ with $k=2$, which means that
the first rational function interpolates $f$ only at the points $-0.9$ and $-0.1$, and the second ones interpolates it at
$+0.1$ and $+1$. These two interpolants and their first and second derivatives agree with that of $f$ at $t=0$.
The solid line in Figure \ref{piecew} corresponds to the curve formed by these two Pad\'e--type rational interpolants.
The two systems have a condition number of $3.25 \times 10^4$ and $2.86 \times 10^4$, respectively.
Then, we construct the
Pad\'{e}--type rational interpolant interpolating $f$ at the 4 points $-0.9,-0.1,+0.1$ and $+1$, and with a
${\cal O}(t^3)$ error at the origin. The system is overdetermined since $l>k$, its condition number is
$1.90 \times 10^3$, and the error is given by the dashed line.
Finally, with the same 4 interpolation points, we construct the interpolant of degree $k=3$. The system is also
overdetermined, its condition number is $6.37 \times 10^{13}$, and we obtain the results given by the dash--dotted line.

\section{Conclusions}

In this paper, we presented in details the particular case of the general rational Hermite interpolation problem
(in rational and barycentric form) where
only the values of the function are interpolated at some points, and where the function and its first derivatives
agree at the origin. Thus, the interpolants constructed in this way possess a Pad\'{e}--type property at 0.
An expression for the error in the real case is given.
The interpolation procedure can be easily
modified to introduce a partial knowledge on the poles and the zeros of the function to approximated.
We also showed how spurious poles can be eliminated. Numerical examples show the interest of the procedures.

\vskip 2mm

The ideas developed in this paper need additional investigations.
An important open problem is to be study
the convergence of the interpolants when the degree tends to infinity as done in \cite{eier} for Pad\'e--type approximants. In our case, we performed some numerical experiments which show that, in some cases, convergence seems to occur while, in some others, no conclusion could be drawn since, for high degrees, the systems are numerically singular.

\vskip 2mm

\noindent {\bf Acknowledgment:} We would like to thank Jean--Paul Berrut for interesting discussions and comments. This work was partially supported by MIUR, PRIN grant no. 20083KLJEZ-003, and by
University of Padova, Project 2010 no. CPDA104492.

\end{document}